\documentclass[11pt,oneside]{amsart}
\usepackage{amsmath}
\usepackage{amsfonts}
\usepackage{amssymb}
\usepackage{amsthm}
\usepackage[colorlinks]{hyperref} 
\usepackage{wasysym}
\usepackage{graphicx}
\usepackage{dashrule}
\usepackage[all]{xy}
\usepackage{multirow}
\usepackage{enumitem}
\usepackage{tikz}
\usetikzlibrary{arrows,decorations.pathmorphing,backgrounds,positioning,fit,calc}
\usepackage{tikz-cd}
\usepackage{framed}

\normalfont\upshape

\theoremstyle{plain}
\newtheorem{defn}{Definition}[section]
\newtheorem{thm}[defn]{Theorem}
\newtheorem{propn}[defn]{Proposition}

\numberwithin{equation}{section}

\theoremstyle{remark}
\newtheorem{rmk}[defn]{Remark}

\DeclareMathOperator{\GL}{GL}
\DeclareMathOperator{\SL}{SL}

\DeclareMathOperator{\spec}{Spec}

\DeclareMathOperator{\conv}{Conv}
\DeclareMathOperator{\stab}{Stab}

\newcommand{\nocontentsline}[3]{}
\newcommand{\tocless}[2]{\bgroup\let\addcontentsline=\nocontentsline#1{#2}\egroup}

\newcommand{\CC}{\mathbb{C}}
\newcommand{\PP}{\mathbb{P}}

\newcommand{\GG}{\mathbb{G}}

\newcommand{\dblslash}{/\! \!/}

\newcommand{\env}{\!
\mathbin{\text{\rotatebox[origin=c]{70}{\scalebox{1.2}{$\approx$}}}} \!}

\newcommand{\inenv}{\dblslash \!_{\circ}}

\newcommand{\ten}{\otimes}

\newcommand{\mf}{\mathfrak}

\newcommand{\kk}{\Bbbk}

\newcommand{\ssfg}{\mathrm{ss,fg}}

\newcommand{\rms}{\mathrm{s}}

\newcommand{\Proj}{{\rm Proj}}
\newcommand{\Spec}{{\rm Spec}}
\newcommand{\nc}{\newcommand}
\nc{\bla}{\phantom{bbbbb}}

\newcommand{\beq}{\begin{equation}}
\newcommand{\eeq}{\end{equation}}
\newcommand{\barr}{\begin{array}}
\newcommand{\earr}{\end{array}}
\newcommand{\beqar}{\begin{eqnarray}}
\newcommand{\eeqar}{\end{eqnarray}}
\newtheorem{theorem}{Theorem}[section]

\newtheorem{definition}[theorem]{Definition}
\newtheorem{remark}[theorem]{Remark}

\newtheorem{exit}[theorem]{Example}

\newenvironment{ex}{\begin{exit}\rm}{\end{exit}}

\nc{\FF}{ {\mathbb F} }
\nc{\HH}{ {\mathbb H} }

\newcommand{\QQ}{{\mathbb Q }}
\newcommand{\UU}{{\mathbb U }}

\newcommand{\calo}{\mathcal{O}}

\newcommand{\calr}{\mathcal{R}}

\nc{\umax}{{U_{\max}}}

\newcommand{\hH}{\hat{H}}
\newcommand{\hU}{\hat{U}}

\nc{\lieq}{{\mathfrak q}}
\nc{\liez}{{\mathfrak z}}
\nc{\lieqs}{{\lieq}^*}
\nc{\lieg}{{\mathfrak g}}
\nc{\liegs}{{\lieg}^*}
\nc{\liep}{{\mathfrak p}}
\nc{\lieps}{{\liep}^*}

\def\a{\alpha}

\def\l{\lambda}

\def\G{\Gamma}

\title{Graded linearisations}

\author{Gergely B\`erczi, Brent Doran,  Frances Kirwan}

\thanks{Early work on this project was supported by the Engineering and Physical Sciences 
Research Council [grant numbers   GR/T016170/1,EP/G000174/1].  Brent Doran was partially supported by Swiss National Science Foundation Award 200021-138071.}

\begin{document}

\begin{abstract}
When the action of a  reductive group on a projective variety has a suitable linearisation, Mumford's geometric invariant theory (GIT) can be used to construct and study an associated quotient variety. In this article we describe how Mumford's GIT can be extended effectively to suitable actions of linear algebraic groups which are not necessarily reductive, with the extra data of a graded linearisation for the action. Any linearisation in the traditional sense for a reductive group action induces a graded linearisation in a natural way. 

The classical examples of moduli spaces which can be constructed using Mumford's GIT are moduli spaces of stable curves and of (semi)stable bundles over a fixed nonsingular curve. This more general construction can be used to construct moduli spaces of unstable objects, such as unstable curves or unstable bundles (with suitable fixed discrete invariants in each case, related to their singularities or Harder--Narasimhan type).

\end{abstract}

\maketitle

In algebraic geometry it is often useful to be able to construct quotients of algebraic 
 varieties by linear algebraic group actions; in particular moduli spaces  (or stacks) 
can be constructed in this way. When the linear algebraic group is reductive, and we have a suitable {linearisation} for its action on a projective variety,  we can use Mumford's geometric invariant theory (GIT) to construct and study such quotient varieties \cite{GIT}. The aim of this article is to describe how Mumford's GIT can be extended effectively to actions of a large family of linear algebraic groups which are not necessarily reductive, with the extra data of a \emph{graded linearisation} for the action. 
Any linearisation in the traditional sense for a reductive group action can be regarded as a graded linearisation in a natural way. 

When a linear algebraic group over an algebraically closed field $\kk$ of characteristic 0 is a semi-direct product $H = U \rtimes R$ of its unipotent radical $U$ and a reductive subgroup $R \cong H/U$ which contains a central one-parameter subgroup $\lambda: \GG_m \to R$ whose adjoint action on the Lie algebra of $U$ has only strictly positive weights, we will see that any linearisation for an action of $H$ on a projective variety $X$ becomes graded if it is twisted by an appropriate (rational) character, and then many of the good properties of Mumford's GIT hold. Many non-reductive linear algebraic group actions arising  in algebraic geometry are actions of groups of this form: for example, any parabolic subgroup of a reductive group has this form, as does the automorphism group of any complete simplicial toric variety \cite{cox}, and  the group of $k$-jets of germs of biholomorphisms of $(\CC^p,0)$ for any positive integers $k$ and $p$ \cite{BK15}.

\begin{ex}
 The automorphism group of the weighted projective plane $\PP(1,1,2)$ 
 with weights 1,1 and 2 is
$$\mbox{Aut}(\PP(1,1,2)) \cong R \ltimes U$$
where $R \cong (GL(2) \times \GG_m)/\GG_m \cong GL(2)$ is reductive and 
$U \cong (\kk^+)^3$ is unipotent
 with elements given by $(x,y,z)  \mapsto (x,y,z+\lambda x^2 + \mu xy + \nu y^2)$ 
for $(\lambda,\mu,\nu) \in \kk^3$.
\end{ex}

\begin{ex}
Under composition modulo $t^{k+1}$ we have a group $\GG_{(k)}$ 
 whose elements are 
$k$-jets of germs of biholomorphisms of $(\CC,0)$:
$$ \{ t \mapsto \phi(t) = a_1 t + a_2 t^2 + \ldots + a_k t^k \,\, \mid \,\, 
a_j \in \CC, a_1 \neq 0 \}.$$
 $\GG_{(k)}$  is isomorphic to a group of matrices of the form
 $$ 
\left\{ \left( \begin{array}{cccc} a_1 & a_2 & \ldots & a_k\\
 0 & (a_1)^2 & \ldots & p_{2k}(\underline{a})\\
 & & \ldots & \\
 0 & 0 & \ldots & (a_1)^k \end{array} \right) : a_1 \in \CC^*, a_2,\ldots a_k \in \CC \right\},$$
where the $(i,j)$th entry $p_{ij}(\underline{a})$ is a polynomial in $a_1, \ldots ,a_k$. 
This reparametrisation group $\GG_{(k)}$ has a one-parameter multiplicative subgroup $\GG_m = \CC^*$ (represented by $\phi(t)=a_1 t$) and  unipotent radical
$\UU_{(k)}$ (represented by $\phi(t)= t + a_2 t^2 + \ldots + a_k t^k$) with
$ \GG_{(k)} \cong \UU_{(k)} \rtimes \CC^*.$ 
\end{ex}

 In Mumford's classical geometric invariant theory the GIT quotient $X/\!/G = \mbox{Proj}({\hat{\calo}}_L(X)^G)$ (where
${\hat{\calo}}_L(X) =   \bigoplus_{k= 0}^{\infty} H^0(X,L^{\otimes k})$)
for an action of a   reductive group $G$ on a  projective variety  $X$ with respect to an ample linearisation $L$ is a projective completion of the geometric quotient $X^s/G$ of the stable set $X^s$. When $X$ is nonsingular then the singularities of $X^s/G$ are very mild, since the stabilisers of stable points are finite subgroups of $G$. If $X^{ss} \neq X^s$ the singularities of $X/\!/G$ are typically more severe, but $X/\!/G$ has a \lq partial desingularisation'  $\tilde{X}/\!/G$ \cite{K2} which is also a projective completion of $X^s/G$ and  is itself  a geometric quotient 
$$\tilde{X}/\!/G = \tilde{X}^{ss}/G$$
by $G$ of an open subset $\tilde{X}^{ss} = \tilde{X}^s$ of a $G$-equivariant blow-up $\tilde{X}$ of $X$. When $X$ is nonsingular then so is $\tilde{X}^{ss}$, and $G$ acts on $\tilde{X}^{ss}$ with finite stabilisers.
$\tilde{X}^{ss}$ is obtained from ${X}^{ss}$ by successively blowing up along the subvarieties of semistable points stabilised by reductive subgroups of $G$ of maximal dimension and then removing the unstable points in the resulting blow-up.

So in the best case in classical GIT we have $X^{ss} = X^s \neq \emptyset$, and then $X^s/G = X/\!/G =  \mbox{Proj}({\hat{\calo}}_L(X)^G)$ is simultaneously a 
 projective variety and a geometric quotient of $X^s$ by the action of $G$. More generally when $X^s \neq \emptyset$ then the geometric quotient $X^s/G$ has a projective completion $\tilde{X}/\!/G$ which is itself a geometric quotient 
$\tilde{X}^{ss}/G$ of an open subset  of a $G$-equivariant blow-up of $X$. Moreover using the 
Hilbert--Mumford criteria for (semi)stability, which allow us to  determine which points of $X$ are stable and which are semistable for the $G$-action without having to know the $G$-invariant sections of powers of $L$, together with the explicit blow-up construction, we can give effective descriptions of $X^s$, $\tilde{X}^{ss} = \tilde{X}^s$ and thus their geometric quotients $X^s/G$ and $\tilde{X}/\!/G$. This is the picture which can be generalised to the action of a non-reductive linear algebraic group, given a graded linearisation of the action.

The immediate problem which arises when trying to extend classical GIT to non-reductive linear algebraic groups $H$ is that in general we cannot define a projective \emph{variety} $X/\!/H = \mbox{Proj}({\hat{\calo}}_L(X)^H)$ 
because ${\hat{\calo}}_L(X)^H$ is not necessarily
finitely generated as a graded algebra, although $\mbox{Proj}({\hat{\calo}}_L(X)^H)$ exists as a scheme. Nonetheless an analogue of classical GIT for non-reductive linear algebraic group actions is described in \cite{BDHK,DK}. Here it is shown that if  $H$ is a linear algebraic group over
$\kk$ acting linearly on a projective variety $X$ with respect to an ample line bundle $L$,
then $X$ has open subvarieties $X^s$
(the locus of \lq stable points') and $X^{ss}$ (\lq semistable points') with a geometric quotient $X^s \to X^s/H$ and an \lq enveloping quotient' $X^{ss} \to X\env H$.
Furthermore  there is a diagram
$$\begin{array}{rcccl}
 & X & - - \rightarrow & \mbox{Proj}({\hat{\calo}}_L(X)^H) &  \\ 
 & \bigcup  & &  \bigcup & \mbox{open} \\
\mbox{semistable} & X^{ss} & \longrightarrow & X\env H & \\
 & \bigcup  & & \bigcup & \mbox{open}\\
\mbox{stable} & X^s & \longrightarrow & X^s/H &
\end{array} $$
where the vertical inclusions are of open subvarieties, and {\em if} ${\hat{\calo}}_L(X)^H$ is finitely generated then $X\env H = \mbox{Proj}({\hat{\calo}}_L(X)^H)$ as in the reductive case. 
However this picture is less helpful than in the case of classical GIT in three significant respects: firstly $X\env H$ is not necessarily a projective variety; secondly (even when ${\hat{\calo}}_L(X)^H$ is finitely generated and so $X\env H = \mbox{Proj}({\hat{\calo}}_L(X)^H)$ is a projective variety) the $H$-invariant morphism 
$X^{ss} \to X\env H $ is {not necessarily a categorical quotient}, and its image is not in general a subvariety of $X \env H$, only a constructible subset; and thirdly there are  in general no obvious analogues of the Hilbert--Mumford criteria for (semi)stability. 

We can see the second of these issues arising in simple examples, when the additive group $\GG_a = \kk^+$ acts on a projective space $\PP^n$ via a linear representation (see Example \ref{simpleex} below). It follows from Jordan canonical form that the representation of $\GG_a$ extends to a representation of $SL(2)$. This enables us to identify $\PP^n \env  \GG_a$ with the reductive GIT quotient
$$(\PP^2 \times \PP^n)/\!/ \SL(2),$$
and thus to see that in general  the quotient morphism 
$q_{\GG_a}: (\PP^n)^{ss,\GG_a} \,\,\,\, {\longrightarrow}  \,\,\,\,\PP^n\env \, \GG_a$ fails to be surjective. 
Twisting the representation of the Borel subgroup $B \cong \GG_a \rtimes \GG_m$ of $SL(2)$   by a character $\chi :B \to \GG_m = \kk^*$ (whose kernel must contain $\GG_a$) 
changes the linearisation but not the action of $B$ on $\PP^n$ to give an enveloping quotient 
$$\PP^n \, \env_\chi \,\, B = (\PP^n \env  \GG_a) /\!/_\chi \,\, \GG_m. $$
It turns out that 
for appropriate choice of (rational) character $\chi$ the complement of the image of $(\PP^n)^{ss, \GG_a}$  in $\PP^n \env \GG_a$ becomes unstable for the $\GG_m$-action and  the morphism 
$$q_B: (\PP^n)^{ss,B,\chi} \,\,\,\, {\longrightarrow}  \,\,\,\,\PP^n\env_\chi \, B$$
to the projective variety $\PP^n\env_\chi \, B$ 
is surjective. This phenomenon occurs more generally.

\begin{definition} \label{defn0.1} Let us call a unipotent linear algebraic group $U$ {\em graded unipotent} 
if there is a homomorphism $\lambda:\GG_m \to Aut(U)$ with the weights of the
 $\GG_m$ action on $Lie(U)$ all {strictly positive}.
For such a homomorphism $\lambda$ let
$$\hat{U} = U \rtimes \GG_m = \{(u,t):u \in U, t \in \GG_m\}$$
be the associated semi-direct product of $U$ and $\GG_m$ with multiplication $(u,t)\cdot(u',t') = (u(\lambda(t)(u')),tt')$. 

When $L$ is very ample, and so induces an embedding of $X$ in a projective space $\PP^n$, we can choose coordinates on $\PP^n$ such that the action of $\GG_m$ on $X$ is diagonal, given by
$$ t \mapsto \left( \begin{array}{cccc} t^{r_0} & 0 & \ldots & 0\\
 0 & t^{r_1} & \ldots & 0\\
 & & \ldots & \\
 0 & 0 & \ldots & t^{r_n} \end{array} \right) $$
where $r_0 \leq r_1 \leq \cdots \leq r_n$. The {\em lowest bounded chamber} for this linear $\GG_m$-action is the closed interval $[r_0,r_j]$ where  $r_0 = \cdots = r_{j-1} <  r_j \leq \cdots \leq r_n$, with interior the open interval $(r_0,r_j)$, unless the action of $\G_m$ on $X$ is trivial; when the action is trivial so that  $r_0 = r_1 = \cdots = r_n$ we will say that $[r_0,r_0]$ is the lowest bounded chamber and it is its own interior. Note that in the situation above, if $\GG_m$ acts trivially then so does $U$.
\end{definition}

Suppose that  $\hat{U}$ acts linearly (with respect to an ample line bundle $L$) on a projective
variety $X$. We can twist  the linearisation  by any (rational) character of $\hat{U}$ (whose kernel must contain $U$) without altering the action of $\hU$ on $X$. In fact by choosing an appropriate rational character we can obtain a GIT picture with many of the good properties of the reductive case, as the following result demonstrates.

\begin{thm}[\cite{BDHK2, BDHK3}] \label{firstthm} Let $U$ be graded unipotent acting linearly on an irreducible projective
variety $X$ with respect to an ample line bundle $L$, and suppose that the linear action extends to $\hat{U} = U \rtimes \GG_m$. Suppose also that semistability coincides with stability in the sense that 
$$ x \in Z_{{\rm min}}  \Rightarrow {\rm Stab}_U(x) = \{ e \} $$
where $Z_{{\rm min}} $ is the union of those connected components of the fixed point set $X^{\GG_m}$ where $\GG_m$ acts on the fibres of $L^*$ with minimum weight. Then the linearisation for the action of $\hat{U}$ on $X$ can be twisted by a rational character of $\hU$ so that 0 lies in the interior of the lowest bounded chamber for the linear $\GG_m$ action on $X$ and\\ 
(i) the algebra ${\hat{\calo}}_{L^{\otimes c}}(X)^{\hat{U}} = \oplus_{m=0}^\infty H^0(X,L^{\otimes cm})^{\hat{U}}$ of $\hat{U}$-invariants is {finitely generated} for any sufficiently divisible integer $c > 0$,
so that the enveloping quotient $X\env \hat{U} = \mbox{Proj}({\hat{\calo}}_{L^{\otimes c}}(X)^{\hat{U}})$ is a projective variety;\\
(ii) $X^{ss,\hat{U}} = X^{s,\hat{U}} $ has a Hilbert--Mumford description as $\bigcap_{u\in U} uX^{s,\GG_m}$, and $X\env \hat{U} = X^{s,\hat{U}} / \hat{U}$ is a geometric quotient of $X^{s,\hat{U}} $ by $\hat{U}$.
 
Moreover, even when the condition that semistability should coincide with stability fails,  
 there is a projective completion of a geometric quotient by $\hU$ of an open subvariety of $X$ (conjecturally $X^{s,\hat{U}} /\hat{U}$), which
 is itself  a geometric quotient  $\tilde{X}^{ss, \hU}/\hat{U}$ by $\hU$ of an open subset $\tilde{X}^{ss, \hU}= \tilde{X}^{s, \hU}$ of a $\hat{U}$-equivariant blow-up $\tilde{X}$ of $X$. 
\end{thm}

If we are interested in constructing quotients of open subsets of $X$ by the action of $U$, not of $\hU$, then we can apply these results  to the diagonal action of $\hat{U} $ on $X \times \PP^1$, where $\hat{U}$ acts on $\PP^1$ via
\begin{equation}  \label{p1action}  (u,t) \mapsto \left(\begin{array}{cc} t & 0 \\ 0 & 1 \end{array}
\right) 
\end{equation}
with kernel $U$, and the linearisation is $L$ tensored with $\mathcal{O}_{\PP^1}(m)$ for $m\! >\!>\! 1$. This gives us a $U$-{invariant open subset} $X^{\hat{s},U}$ of $X$ 
 with a {geometric quotient} $X^{\hat{s},U}/U$ by $U$ which is isomorphic to the geometric quotient by $\hU$ of the open subvariety $\GG_m(X^{\hat{s},U} \times \{ [1:1]\} $ of $X \times \PP^1$; moreover it  has a {projective completion} which
 is a geometric quotient by $\hat{U}$ of an open subvariety of a $\hat{U}$-equivariant blow-up of $X \times  \PP^1$. 
Furthermore in this set-up there are
Hilbert--Mumford-like criteria for (semi)stability.

This motivates the following definitions.

 \begin{definition} \label{defnextendedlin} An {\em extended linearisation} $\mathcal{L}$ of an action of a linear algebraic group $H$ on a projective variety $X$ is given by the data:\\
(a) a line bundle $L$ on $X$;\\
(b) a semi-direct product $\hat{H} = H \rtimes \GG_m$ of $H$ by $\GG_m$; \\
(c) an extension of the $H$-action on $X$ to $\hat{H}$ and a lift of the $\hat{H}$-action to $L$.

Given an extended linearisation $\mathcal{L}$ and a rational number $q \in \QQ$, define the \lq $q$-hat-stable' locus $X^{\hat{s},q} = X^{\hat{s},q,\mathcal{L}}$ to be the $H$-{invariant open subvariety}  of $X$ determined by
$$  X^{\hat{s},q} \times  \{ [1:1]\}  = (X \times \PP^1)^{s,\hat{H}} \cap (X \times \{ [1:1]\} )$$
where $\hat{H}$ acts on $\PP^1$ as at (\ref{p1action}) above with its linearisation on
$\mathcal{O}_{\PP^1}(1)$  twisted by $q$, and
$(X \times \PP^1)^{s,\hat{H}}$ is defined 
with respect to the induced linearisation for the $\hH$-action on $X \times \PP^1$ on $L$ tensored with $\mathcal{O}_{\PP^1}(m)$ for $m\! >\!>\! 1$. 

We then have a {geometric quotient} $X^{\hat{s},q}/H$ by $H$ which is isomorphic to an open subvariety of $(X \times \PP^1)^{s,\hat{H}} / \hat{H}$. 
\end{definition}

\begin{rmk}
Given a linearisation in the classical sense of an action of a linear algebraic group $H$ on a projective variety $X$ with respect to a line bundle $L$, we have a \lq trivial extended linearisation' for which $\hat{H} = H \times \GG_m$ and $\GG_m$ acts trivially on $X$ and on $L$. Then if $q \in (0,1)$ the $q$-hat-stable locus $X^{\hat{s},q}$ coincides with the stable locus defined as in \cite{BDHK} for the action of $H$ on $X$ with the given linearisation, while if $q \notin [0,1]$ the $q$-hat-stable locus $X^{\hat{s},q}$ is empty.
\end{rmk}

\begin{definition} \label{defngradedlin}
A {\em  linear algebraic group with graded unipotent radical} is a linear algebraic group $H$ with unipotent radical $U$, equipped with a semi-direct product $\hH = H \rtimes \GG_m$ such that the adjoint action of $\GG_m$ on the Lie algebra of $U$ has only strictly positive weights and the induced conjugation action of $\GG_m$ on $H/U$ is trivial. 

A {\em graded linearisation} $\mathcal{L}$ of an action of  $H$ on a projective variety $X$ is then an extended linearisation in the sense of Definition \ref{defnextendedlin} for this choice of $\hH$, such that the $\hH$-linearisation is twisted  by a rational character of $\hH$ so that 0 lies in the interior of the lowest bounded chamber for the $\GG_m$ action;  we will assume that the line bundle $L$ on $X$ is ample unless stated otherwise. Given a graded linearisation $\mathcal{L}$,  the \lq hat-stable' locus $X^{\hat{s}} = X^{\hat{s},\mathcal{L}}$ is the $0$-hat-stable locus $X^{\hat{s},0}$ as defined in Definition \ref{defnextendedlin} when $q=0$.
\end{definition}

\begin{rmk}
When $H$ is a linear algebraic group with graded unipotent radical $U$ and $\mathcal{L}$ is a graded linearisation for an action of $H$ on a projective variety $X$ (with respect to an ample line bundle $L$ on $X$), then we can apply Theorem \ref{firstthm} to the action of $\hU$ on $X \times \PP^1$ as above, and then apply classical GIT and the partial desingularisation construction of \cite{K2} to the induced action of the reductive group $\hH/\hU \cong H/U$. Thus the {geometric quotient} $X^{\hat{s}}/H$ by $H$ has a {projective completion} which
 is a geometric quotient by $\hat{H}$ of an open subset of a $\hat{H}$-equivariant blow-up of $X \times  \PP^1$. 
Furthermore  the {geometric quotient} $X^{\hat{s}}/H$ by $H$ and its {projective completion} can be described using
Hilbert--Mumford-like criteria combined with an explicit blow-up construction. 
\end{rmk}

\begin{rmk}
Definitions \ref{defnextendedlin}  and \ref{defngradedlin} can be extended to define $T$-extended linearisations and $T$-graded linearisations for the actions of linear algebraic groups with $T$-graded unipotent radical, for any torus $T$.
\end{rmk}

The layout of this  article is as follows. In $\S$1 we will review GIT with classical linearisations \cite{BDHK,DK,GIT}. In $\S$2 we will describe extended, graded and torus-graded linearisations and the associated geometric invariant theory for these. Finally $\S$3 describes some potential applications, including 
the construction of moduli spaces of \lq unstable' objects, such as unstable projective curves or unstable sheaves over a fixed nonsingular projective variety (with suitable fixed discrete invariants in each case, involving their singularities or Harder--Narasimhan type).

\section{GIT with classical linearisations}

\subsection{GIT for reductive groups}

In Mumford's classical {G}eometric {I}nvariant {T}heory a linearisation (more precisely, an ample linearisation) of an action of a   reductive group $G$ on an irreducible   projective variety  $X$ over an algebraically closed field $\kk$ of characteristic 0 
is given by  an ample line bundle $L$ on $X$ and a lift of the action to $L$;
when $X$ is embedded in a projective space $\PP^n$ and $L = \calo(1)$,   the action is given by a
representation $\rho:G \to GL(n+1)$ and $ {\hat{\calo}}_L(X) =  \bigoplus_{k= 0}^{\infty} H^0(X,L^{\otimes k})$ is $\kk[x_0,\ldots,x_n]/\mathcal{I}_X$ where $\mathcal{I}_X$ is the ideal generated by the homogeneous polynomials which vanish on $X$. We consider the picture:
$$\begin{array}{ccccl}
(X,L) & \leadsto 
 & {\hat{\calo}}_L(X)&=&  \bigoplus_{k= 0}^{\infty} H^0(X,L^{\otimes k})\\
| &&&& \\
| & & \bigcup \!| & &\\
\downarrow & & & & \\
X/\!/G & {\reflectbox{\ensuremath{\leadsto}}}   
 & {\hat{\calo}}_L(X)^G & & \mbox{ algebra of invariants. }
\end{array}$$
Since $G$ is reductive, the algebra of $G$-invariants ${\hat{\calo}}_L(X)^G$ is  {\em finitely generated} as a graded  algebra with associated projective variety  
$X/\!/G = \mbox{Proj}({\hat{\calo}}_L(X)^G)$. The inclusion of ${\hat{\calo}}_L(X)^G$ in ${\hat{\calo}}_L(X)$ determines a rational map $X - - \rightarrow X/\!/G$ which fits into a diagram
$$\begin{array}{rcccl}
 & X & - - \rightarrow & X/\!/G & \mbox{ projective variety}\\
 & \bigcup  & & || & \\
\mbox{semistable} & X^{ss} & \stackrel{\mbox{onto}}{\longrightarrow} & X/\!/G & \\
 & \bigcup  & & \bigcup & \mbox{open}\\
\mbox{stable} & X^s & \longrightarrow & X^s/G &
\end{array} $$     
where $X^s$ and $X^{ss}$ are open subvarieties of $X$, the GIT quotient ${X}/\!/G$ is a categorical quotient for the action of $G$ on $X^{ss}$ via the  $G$-invariant surjective morphism $\phi_G: X^{ss} \to X/\!/G$, and  
$$\phi_G(x) = \phi_G( y) \Leftrightarrow \overline{Gx} \cap \overline{Gy} \cap X^{ss} \neq \emptyset.$$
\begin{rmk}
When $\kk = \CC$ then $G$ is reductive if and only if it is the complexification $G=K_\CC$ of a maximal compact subgroup $K$, 
 and then $X/\!/G = \mu^{-1}(0)/K$ for a suitable \lq moment map' $\mu$ for the action of $K$.
\end{rmk}

The subsets $X^{ss}$ and $X^s$ of $X$ for a linear action of a reductive group $G$ with respect to an ample linearisation are characterised by the following properties (see
\cite[Chapter 2]{GIT}, \cite{New}). 

\begin{propn} (Hilbert--Mumford criteria for reductive group actions)\\
\label{sss} (i) A point $x \in X$ is semistable (respectively
stable) for the action of $G$ on $X$ if and only if for every
$g\in G$ the point $gx$ is semistable (respectively
stable) for the action of a fixed maximal torus $T$ of $G$.

\noindent (ii) A point $x \in X$ with homogeneous coordinates $[x_0:\ldots:x_n]$
in some coordinate system on $\PP^n$
is semistable (respectively stable) for the action of a maximal 
torus $T$ of $G$ acting diagonally on $\PP^n$ with
weights $\a_0, \ldots, \a_n$ if and only if the convex hull
$$\conv \{\a_i :x_i \neq 0\}$$
contains $0$ (respectively contains $0$ in its interior).
\end{propn}

 The GIT quotient $X/\!/G$ is a projective completion of the geometric quotient $X^s/G$ of the stable set $X^s$. When $X$ is nonsingular then the singularities of $X^s/G$ are very mild, since the stabilisers of stable points are finite subgroups of $G$. If $X^{ss} \neq X^s \neq \emptyset$ the singularities of $X/\!/G$ are typically more severe, but $X/\!/G$ has a \lq partial desingularisation'  $\tilde{X}/\!/G$ which is also a projective completion of $X^s/G$ and  is itself  a geometric quotient 
$$\tilde{X}/\!/G = \tilde{X}^{ss}/G$$
by $G$ of an open subset $\tilde{X}^{ss} = \tilde{X}^s$ of a $G$-equivariant blow-up $\tilde{X}$ of $X$ \cite{K2}. 

$\tilde{X}^{ss}$ is obtained from ${X}^{ss}$ by successively blowing up along the subvarieties of semistable points stabilised by reductive subgroups of $G$ of maximal dimension and then removing the unstable points in the resulting blow-up, as follows. 
We suppose  that $X$ has some stable points. There exist 
semistable points of $X$ which are not stable if and only if there exists a 
non-trivial connected reductive subgroup of $G$ fixing a semistable point. Let 
$r>0$ be the maximal dimension of a reductive subgroup of $G$ 
fixing a point of $X^{ss}$ and let $\calr(r)$ be a set of representatives of conjugacy 
classes of all connected reductive subgroups $R$ of 
dimension $r$ in $G$ such that 
$$ Z^{ss}_{R} = \{ x \in X^{ss} :  \mbox{$R$ fixes $x$}\} $$
is non-empty. Then
$$
\bigcup_{R \in \calr(r)} GZ^{ss}_{R}
$$
is a disjoint union of nonsingular closed subvarieties of $X^{ss}$. The action of 
$G$ on $X^{ss}$ lifts to an action on the blow-up $X_{(1)}$ of 
$X^{ss}$ along $\bigcup_{R \in \calr(r)} GZ_R^{ss}$ which can be linearised so that the complement 
of $X_{(1)}^{ss}$ in $X_{(1)}$ is the proper transform of the 
subset $\phi^{-1}(\phi(GZ_R^{ss}))$ of $X^{ss}$ where $\phi:X^{ss} \to X/\!/G$ is the quotient 
map (see \cite{K2} 7.17). 
Here we use the linearisation with respect to (a tensor power of) the pullback of the ample line bundle $L$ on $X$ perturbed by a sufficiently small multiple of the exceptional divisor $E_{(1)}$. This will give us an ample line bundle on the blow-up  $\psi:X_{(1)} \to X$ , and if the perturbation is sufficiently small it will have the property that 
$$   \psi^{-1}(X^{s}) \subseteq  X_{(1)}^{s} \subseteq  X_{(1)}^{ss} \subseteq \psi^{-1}(X^{ss}) = X_{(1)},$$
and the stable and semistable subsets $X_{(1)}^{s}$ and $X_{(1)}^{ss}$ will be independent of the choice of perturbation.
Moreover no point of $X_{(1)}^{ss}$ is fixed by a 
reductive subgroup of $G$ of dimension at least $r$, and a point in $ X_{(1)}^{ss} $ 
is fixed by a reductive subgroup $R$ of 
dimension less than $r$ in $G$ if and only if it belongs to the proper transform of the 
subvariety $Z_R^{ss}$ of $X^{ss}$.

\begin{rmk} \label{cfK2} In fact in \cite{K2} $X$ itself is blown  up along the closure $\overline{\bigcup_{R \in \calr(r)} GZ^{ss}_{R} }$ of $\bigcup_{R \in \calr(r)} GZ^{ss}_{R}$ in $X$ (or in a projective completion of $X^{ss}$ with a $G$-equivariant morphism to $X$ which is an isomorphism over $X^{ss}$). This gives us  a projective variety $\bar{X}_{(1)}$ and blow-down map $\bar{\psi}: \bar{X}_{(1)} \to X$ restricting to $\psi:X_{(1)} \to X$ where $\bar{\psi}^{-1}(X^{ss}) = X_{(1)}$. We can then choose a sufficiently small perturbation of the pullback to $\bar{X}_{(1)}$ of the linearisation on $X$ which provides an ample linearisation of the projective variety $\bar{X}_{(1)}$ such that
$   \bar{\psi}^{-1}(X^{s}) \subseteq \bar{X}_{(1)}^{s} \subseteq \bar{X}_{(1)}^{ss} \subseteq \bar{\psi}^{-1}(X^{ss}) = X_{(1)},$
and moreover the restriction of the linearisation to $X_{(1)}$ is obtained from the pullback of $L$ by perturbing  by a sufficiently small multiple of the exceptional divisor $E_{(1)}$.
\end{rmk}

If $r>1$ the same procedure can be applied to $X_{(1)}^{ss}$ to obtain $X_{(2)}^{ss}$ such 
that no reductive subgroup of $G$ of dimension at least 
$r-1$ fixes a point of $X_{(2)}^{ss}$. After repeating this enough times, we obtain 
$X_{(0)}^{ss} = X^{ss},X_{(1)}^{ss},X_{(2)}^{ss},\ldots,X_{(r)}^{ss}=\tilde{X}^{ss}$ such that
no reductive subgroup of $G$ of positive dimension fixes a point of $\tilde{X}^{ss}$.  Similarly 
$\tilde{X}/\!/G = \tilde{X}^{ss}/G$ can be obtained from $X/\!/G$ by blowing 
up along the proper transforms of the images $Z_R /\!/N$ 
in $X/\!/G$ of the subvarieties $GZ_R^{ss}$ of $X^{ss}$ in decreasing order of $\dim R$.

Thus when a reductive group $G$ acts linearly on an irreducible projective variety $X$ with respect to an ample linearisation, we can summarise the GIT output when $X^s \neq \emptyset$ as follows:\\
i) The best case is when $X^{ss} = X^s \neq \emptyset$, and then the GIT quotient $X/\!/G = X^s/G$ is a projective variety which is a geometric quotient of the open subvariety $X^s$ of $X$. Furthermore  the stabiliser in $G$ of every $x \in X^s$ is finite, so if $X$ is nonsingular then $X/\!/G$ has at worst orbifold singularities.\\
ii) When  $X^{ss} \neq X^s \neq \emptyset$ then the GIT quotient $X/\!/G$ is a projective completion of the geometric quotient $X^s/G$. Typically the singularities of $X/\!/G$ are significantly more serious than those of $X^s/G$, but $X^s/G$ has another projective completion $\tilde{X}/\!/G = \tilde{X}^{s}/G$ which is a \lq partial desingularisation' of $X/\!/G$ in the sense described above.

\subsection{GIT for non-reductive groups}

Now let $X$ be an irreducible projective variety over an algebraically closed field $\kk$ of characteristic 0 and let $H$ be a  linear algebraic group, with unipotent radical $U$, acting on $X$ with  an ample  {linearisation} of the action (that is, an ample
line bundle $L$ on $X$ and a lift 
 of the action to $L$). First we will define stability and semistability for the linear action of the unipotent group $U$.

\begin{definition} (cf. \cite{DK} $\S$4 and \cite{DK} 5.3.7). \label{defnssetc}
Let $I = \bigcup_{m>0} H^0(X,L^{\otimes m})^U$
and for $f \in I$ let $X_f$ be the $U$-invariant affine open subset
of $X$ where $f$ does not vanish, with ${\calo}(X_f)$ its coordinate ring. 
A point $x \in X$ is called {\em semistable} for the linear action of the unipotent group $U$  if 
 there exists some $f \in I$
which does not vanish at $x$, and such that $\calo(X_f)^U$ is finitely generated as a  graded algebra.
The $U$-{\em semistable locus} of $X$ is  $X^{ss,
U} =  \bigcup_{f \in I^{fg}} X_f$ where
$$I^{fg} = \{f
\in I \ | \ {\calo}(X_f)^U
\mbox{ is finitely generated }   \}.$$
The  {\em stable locus}
 of $X$ for the linear action of $U$  is
      $ X^{s,U} =
\bigcup_{f \in I^{lts} } X_f$ where
\begin{eqnarray*} I^{lts}\ \  =\ \  \{f
\in I^{fg}  \ | \  
  \mbox{the quotient map }  q_U: X_f \longrightarrow \Spec({\calo}(X_f)^U) \mbox{ is a locally trivial
geometric quotient} \}. \end{eqnarray*}
\label{defn:envelopquot}
The {\em enveloped quotient} of
$X^{ss,U}$ by the linear $U$-action is $q_U: X^{ss, U} \rightarrow q_U(X^{ss,U})$, where
$q_U: X^{ss, U} \rightarrow \Proj({\hat{\calo}}_L(X)^U)$ is the natural
morphism of schemes and
$q_U(X^{ss,U})$ is a dense constructible subset of the {\em
enveloping quotient}
$$X \env  U = \bigcup_{f \in I^{ss,fg}}
\Spec({\calo}(X_f)^U)$$ of $X^{ss, U}$. 
\end{definition}

\begin{rmk} 
If ${\hat{\calo}}_L(X)^U$ is finitely generated then $X\env U$ is the projective
variety $\Proj({\hat{\calo}}_L(X)^U)$. Note that even in this case $q_U(X^{ss,U})$ is 
not necessarily a subvariety of $X \env U $ (see for  example \cite{DK} $\S$6).

   \label{remclaim} The enveloping quotient 
$X \env  U$ has quasi-projective open subvarieties (\lq inner enveloping quotients' $X \inenv U$) which contain the enveloped quotient $q_U(X^{ss})$  and have ample line bundles pulling back to  positive tensor powers of $L$
under the natural map $q_U:X^{ss} \to X\env U$ (see \cite{BDHK} for details). 
\end{rmk}

The $H$-semistable set
 $X^{ss} = X^{ss,H}$, enveloped and enveloping quotients and inner enveloping quotients
$$q_H: X^{ss} \to q_H(X^{ss}) \subseteq X \inenv H \subseteq X \env H$$
for the linear action of $H$ are defined exactly  as for the unipotent case in Definition \ref{defnssetc} and Remark \ref{remclaim} (cf. \cite{BDHK}). However the definition given in \cite{BDHK} of the stable set $X^s= X^{s,H}$ for the linear action of $H$ combines the unipotent and reductive cases as follows.

\begin{definition} \label{def:GiSt1.1}
Let $H$ be a linear algebraic group acting on an irreducible variety $X$ and $L \to X$ a
linearisation for the action. The \emph{stable locus} is the open subvariety
\[
X^{\rms}= \bigcup_{f \in I^{\rms}} X_f
\]
of $X^{ss}$, where $I^{\rms} \subseteq \bigcup_{r>0} H^0(X,L^{\ten r})^H$ is the subset
of $H$-invariant sections $f$ of tensor powers of $L$ satisfying the following conditions:
\begin{enumerate}
\item \label{itm:GiSt1.1-1} the open set $X_f$ is affine (this is automatically true when $X$ is projective);
\item \label{itm:GiSt1.1-2} the action of $H$ on $X_f$ is closed with all
  stabilisers finite groups; and
\item \label{itm:GiSt1.1-3} the restriction of the $U$-enveloping quotient map
 \[
q_{U}:X_f \to \spec((S^{U})_{(f)})
\]
is a principal $U$-bundle for the action of $U$ on $X_f$. 
\end{enumerate}
\end{definition}

\begin{rmk}
When $H$ is reductive or unipotent these definitions of $X^{ss,H}$ and $X^{s,H}$ coincide with those already given.
\end{rmk}

\begin{ex} \label{simpleex} Let $\GG_a= \kk^+$ act linearly  on $\PP^n$ via a representation on $\kk^{n+1}$.
We can choose coordinates in which the generator of $\mathrm{Lie}(\GG_a)$
has { Jordan normal form} with blocks of size $k_1 +1,\ldots,k_q+1$.
The linear $\GG_a$ action therefore extends to the reductive group $G = SL(2)$ with
$$\GG_a = \left\{ \left(\begin{array}{cc} 1& a \\ 0 & 1 \end{array}
\right) : a \in \kk \right\} \leqslant G$$
via the identification $\kk^{n+1} \cong \bigoplus_{i=1}^q Sym^{k_i}(\kk^2)$.
In fact in this case the $\GG_a$-invariants are finitely generated by the Weitzenb\"ock theorem \cite{Dolg}, so we have 
$$\PP^n \env \GG_a = \mbox{{Proj}}((\kk[x_0,\ldots,x_n])^{\GG_a}).$$
The Weitzenb\"ock theorem can be proved by considering the identification of $G$-spaces 
$$G \times_{\GG_a} \PP^n \cong (G/\GG_a) \times \PP^n \cong (\kk^2 \setminus \{0\}) \times \PP^n$$ via $(g,x) \mapsto (g\GG_a,gx)$,
composed with the inclusions $ (\kk^2 \setminus \{0\}) \times \PP^n 
\subseteq \kk^2 \times \PP^n \subseteq \PP^2 \times \PP^n$.
We choose a linearisation for the diagonal $G$-action on $\PP^2 \times \PP^n$ given by  $L = \mathcal{O}_{\PP^n}(1)$ tensored with $\mathcal{O}_{\PP^1}(m)$ for $m\! >\!>\! 1$. Then restricting $G$-invariant sections of tensor powers of this linearisation to $\{1\} \times \PP^n$ defines an isomorphism onto the algebra of $\GG_a$-invariant sections of tensor powers of $L$, and we have  
$$\PP^n \env \GG_a  = \mbox{{Proj}}((\kk[x_0,\ldots,x_n])^{\GG_a} \cong (\PP^2 \times \PP^n)/\!/SL(2).$$
We can see how to interpret $ (\PP^n)^{s} =  (\PP^n)^{s,\GG_a}$ and $ (\PP^n)^{ss} = (\PP^n)^{ss,\GG_a}$ as well as 
 the morphism $ (\PP^n)^{ss,\GG_a}  \longrightarrow \PP^n \env \GG_a$ from the diagram 
$$ 
{\begin{array}{cccc}
\PP^2 \times \PP^n & - - \rightarrow & \PP^2 \times \PP^n/\!/G \\
  \bigcup & &  || \\
   \PP^n \cong \{[1:0:1]\} \times \PP^n  & - - \to &  \PP^n \env \GG_a\\
 \bigcup   & &  || \\
 (\PP^n)^{ss}  & \longrightarrow  & \PP^n\env \GG_a & \\
\bigcup  & &  \bigcup \\
 (\PP^n)^s & \longrightarrow & (\PP^n)^s/\GG_a. &
\end{array}}$$  
In particular  the morphism $ (\PP^n)^{ss} 
 \longrightarrow \PP^n \env \GG_a$
 is {\em not} onto when
$$\PP^n = \PP(Sym^n(\kk^2)) = \{ \mbox{ n unordered points on $\PP^1$}\}$$  for $n \geq 3$. When $n=3$
then $(\PP^3)^{ss} = (\PP^3)^s = 
 \{ \mbox{ 3 unordered points on $\PP^1$, at most one at $\infty$}\}$
while its image in
$\PP^3\env \GG_a = (\PP^3)^s/\GG_a \,\, \sqcup \,\, \PP^3/\!/SL(2)$
is the open subset $(\PP^3)^s/\GG_a$ which does not include the
\lq boundary' points coming from $0 \in \kk^2 \subseteq \PP^2.$ When $n=4$ then  $(\PP^4)^{ss} \neq (\PP^4)^s$ and the image of $(\PP^4)^{ss}$ in $\PP^4 \env \GG_a$ is a constructible subset but not a subvariety.

Let $$ 
B = \left\{  \left(\begin{array}{cc} a& b \\ 0 & a^{-1} \end{array}
\right) : a \in \GG_m, b \in \kk \right\} \cong \GG_a \rtimes \GG_m$$
be the standard Borel subgroup of $SL(2)$, acting on $\PP^n$ via a linear representation on $\kk^{n+1}$. 
Then $\kk^{n+1} \cong \bigoplus_{i=1}^q Sym^{k_i}(\kk^2) \otimes \kk_{(r_i)}$
where $B$ acts on $\kk_{(r)} = \kk$ as multiplication by a character
$\chi_r.$
Twisting the representation of $B$ on $\kk^{n+1}$ by a character $\chi$ 
changes the linearisation but not the action of $B$ on $\PP^n$ to give
$$\PP^n \env_\chi B = (\PP^n \env  \GG_a) /\!/_\chi \GG_m. $$
For appropriate $\chi$, in the example above the \lq boundary points' in $\PP^3 \env \GG_a$ become unstable for the $\GG_m$ action and  we have a surjective morphism 
$$ (\PP^3)^{ss,B,\chi} \longrightarrow\PP^3\env_\chi B.$$
It turns out, as will be discussed next, that this is a special case of a more general phenomenon.
\end{ex}

\subsection{GIT for linear algebraic groups with graded unipotent radicals}

Recall from Definition \ref{defngradedlin} that a linear algebraic group with graded unipotent radical is a linear algebraic group $H$ with unipotent radical $U$, equipped with a semi-direct product $\hH = H \rtimes \GG_m$ such that the adjoint action of $\GG_m$ on the Lie algebra of $U$ has only strictly positive weights and the induced conjugation action of $\GG_m$ on $H/U$ is trivial. 

\begin{rmk}
Suppose that $H = U \rtimes R$ where the reductive group $R = H/U$ itself contains a central one-parameter subgroup whose conjugation action on the Lie algebra of $U$ has all weights strictly positive. Then   corresponding semi-direct products $\hat{U}$and $\hH$ can be constructed such that $\hU$ is isomorphic to a subgroup of $H$, and any linear action of $H$ on a projective variety $X$ can be extended to a linear action of $\hH$. We will call this situation an \lq internal grading' for the unipotent radical of $H$. We will call this situation an \lq internal grading' for the unipotent radical of $H$.
\end{rmk}

Given any action of ${\hat{H}}$ on a projective variety $X$ which is linear with respect to an ample line bundle $L$ on $X$, it is shown in \cite{BDHK2, BDHK3} that {provided}:

(i)  we are willing to replace $L$ with a suitable tensor power $L^{\otimes m}$, with $m\geq 1$ sufficiently divisible, and to twist the linearisation of the action of ${\hat{H}}$ by a suitable (rational) character of ${\hat{H}}$ with kernel containing $H$, and moreover

(ii)  \lq semistability coincides with stability' for the action of the unipotent radical ${U}$,

\noindent  then
 the ${\hat{H}}$-invariants form a finitely generated algebra. Moreover in this situation the natural quotient morphism $q_H$ from the semistable locus $X^{ss,{\hat{H}}}$ to the enveloping quotient $X\env {\hat{H}}$ is surjective, and  expresses the projective variety $X\env {\hat{H}}$ as a categorical quotient of $X^{ss,{\hat{H}}}$. Furthermore this locus $X^{ss,{\hat{H}}} = X^{s,{\hat{H}}}$  can be described using Hilbert--Mumford criteria.

In  \cite{BDHK3} it is also shown that when the condition that semistability coincides with stability for the unipotent radical is not satisfied, but is replaced with the weaker condition that the stabiliser in $U$ of a generic point in $X$ is trivial, or equivalently
\begin{equation}  
\text{$\min_{x \in X} \, \dim(\stab_{U}(x)) = 0,$} 
\end{equation}
then
there is a  sequence of blow-ups of $X$ along $\hat{H}$-invariant subvarieties  (similar to that of \cite{K2} when $H$ is reductive) resulting in a projective variety $\hat{X}$ with an induced linear action of $\hat{H}$ satisfying the  condition that semistability coincides with stability for the unipotent radical $U$. In this way we obtain a projective variety $\widehat{X \times \PP^1} \env \hat{H}$ which is a categorical quotient by $\hat{H}$ of a $\hat{H}$-invariant open subset of a blow-up of ${X} \times \kk$ and contains as an open subset a geometric quotient of an $H$-invariant open subset $X^{\hat{s},H}$ of $X$ by $H$, where the geometric quotient $X^{\hat{s},H}/H$ and the projective variety $\widehat{X \times \PP^1} \env \hat{H}$ have descriptions in terms of Hilbert--Mumford-like criteria and  the explicit blow-up construction. 

\begin{rmk}
In fact this can be generalised to the case when $\min_{x \in X} \,\, \dim(\stab_{U}(x))>0$  \cite{BDHK3,behjk16, bejk16}.
\end{rmk}

The description of the condition we need the action of the unipotent radical $U$ of $H$ to satisfy as  \lq semistability coincides with stability' is a rather loose one. To describe it  more precisely,
let $L \to X$ be a very ample linearisation
of  the action of $\hat{H}$ on an irreducible projective variety $X$. 
Let $\chi: \hat{H} \to \GG_m$ be a character of $\hat{H}$ with kernel containing $H$;  such characters $\chi $ can be identified with integers so that the integer 1 corresponds to the character which fits into the exact sequence $H \to \hat{H} \to \GG_m$. Let $\omega_{\min}$ be the minimal weight for the $\GG_m$-action on
$V:=H^0(X,L)^*$ and let $V_{\min}$ be the weight space of weight $\omega_{\min}$ in
$V$. Suppose that $\omega_{\min} < \omega_{\min +1} < 
\cdots < \omega_{\max} $ are the weights with which the one-parameter subgroup $\GG_m \leq {\hat{U}} \leq \hat{H}$ acts on the fibres of the tautological line bundle $\mathcal{O}_{\PP((H^0(X,L)^*)}(-1)$ over points of the connected components of the fixed point set $\PP((H^0(X,L)^*)^{\GG_m}$ for the action of $\GG_m$ on $\PP((H^0(X,L)^*)$; since $L$ is very ample $X$ embeds in $\PP((H^0(X,L)^*)$ and the line bundle $L$ extends to the dual $\mathcal{O}_{\PP((H^0(X,L)^*)}(1)$ of the tautological line bundle $\mathcal{O}_{\PP((H^0(X,L)^*)}(-1)$.
Without loss of generality we may assume that there exist at least two distinct such weights, since otherwise the action of the unipotent radical $U$ of $H$ on $X$ is trivial, and so the action of $H$ is via an action of the reductive group $R=H/U$ and reductive GIT can be applied.

Let $\chi$ be a rational character of $\GG_m$ (lifting to a rational character of $\hH$ as above)  such that 
$$ \omega_{\min} <
\chi < \omega_{\min + 1}; $$
 we will call rational characters $\chi$  with this property {\em  adapted} to the linear action of $\hat{H}$, and we will call the linearisation adapted if $\omega_{\min} <0 < \omega_{\min + 1} $;  we will call $\chi$  {\em borderline adapted } to the linear action of $\hat{H}$ if $\chi = \omega_{\min}$, and the linearisation borderline adapted if $\omega_{\min} = 0$. 
 The linearisation of the action of $\hat{H}$ on $X$ with respect to an ample line bundle $L^{\otimes c}$ for a sufficiently divisible positive integer $c$ such that $c\chi$ is a character can be twisted by this character; effectively the weights $\omega_j$ are replaced with $\omega_j-\chi$ and this twisted linearisation is adapted in the sense above. 
Let $X^{s,\GG_m}_{\min+}$ denote the stable set in $X$ for the linear action of $\GG_m$ with respect to this adapted linearisation  and for a maximal torus $T$ of $\hH$ containing $\GG_m$, let $X^{s,T}_{\min+}$ denote the stable set in $X$ for the linear action of $T$ with respect to the adapted linearisation; by the theory of variation of (classical) GIT \cite{dh98,Thaddeus},   $X^{s,\GG_m}_{\min+}$ and $X^{s,T}_{\min+}$ are independent of the choice of adapted rational character $\chi$.
Let $$X^{s,{\hat{U}}}_{\min+} = X \setminus {\hat{U}} (X \setminus X^{s,\GG_m}_{\min+}) = \bigcap_{u \in U} u X^{s,\GG_m}_{\min+}$$ be the complement of the ${\hat{U}}$-sweep (or equivalently the $U$-sweep) of the complement of $X^{s,\GG_m}_{\min+}$,
and let
$$X^{s,{\hat{H}}}_{\min+} =  \bigcap_{h \in H} u X^{s,T}_{\min+},$$
while
\[
Z_{\min} =X \cap \PP(V_{\min})=\left\{
\begin{array}{c|c}
\multirow{2}{*}{$x \in X$} & \text{$x$ is a $\GG_m$-fixed point and} \\ 
 & \text{$\GG_m$ acts on $L^*|_x$ with weight $\omega_{\min}$} 
\end{array}
\right\}
\]
and
\[
X^0_{\min} =\{x\in X \mid \lim_{t \to 0, \,\, t \in \GG_m} t \cdot x \in Z_{\min}\}.
\]   
Note that $X^0_{\min}$ is $\hat{U}$-invariant and $X^{s,{\hat{U}}}_{\min+}  = X^0_{\min} \setminus U Z_{\min}$.

The condition that \lq semistability coincides with stability' for the linear action of ${\hat{U}}$ required in \cite{BDHK2} is slightly stronger than that required in \cite{BDHK3}; in \cite{BDHK3}  the hypothesis needed for the
$\hat{U}$-linearisation $L \to X$ is that 

\begin{equation}  
\text{$\stab_{U}(z) = \{ e \} $ for every $z \in Z_{\min}$}. \tag{$\mf{C}^*$}
\end{equation}

\begin{thm} \label{mainthm}  {\rm \cite{BDHK3}}
Let $H$ be a linear algebraic group over $\kk$ with unipotent radical $U$.
 Let $\hat{H} = H \rtimes \GG_m$ be a semidirect product of $H$ by $\GG_m$ with subgroup $\hat{U} = U \rtimes \GG_m$,
where the conjugation action of $\GG_m$ on $U$ is such that all the weights
of the induced $\GG_m$-action on the Lie algebra of $U$ are strictly positive, while the induced conjugation action of $\GG_m$ on $R=H/U$ is trivial.
Suppose that ${\hat{H}}$  acts linearly on an irreducible projective variety $X$ with respect to an  ample line bundle $L$, and that 
the linearisation is adapted in the sense above. 
Suppose also that the linear action of ${\hat{U}}$ on $X$ satisfies the condition  ($\mf{C}^*$). 
Then 
    
(i) the open subvariety $X^{s,{\hat{U}}}_{\min+}$ of $X$ has a geometric quotient $X \env \hU = X^{s,{\hat{U}}}_{\min+}/\hat{U}$ by $\hU$ which is a projective variety, while
 
(ii) the open subvariety $X^{s,{\hat{H}}}_{\min+}$ of $X$ has a categorical quotient $X \env \hH$ by $\hH$ which is also a projective variety.
\end{thm} 

\begin{rmk} In order to prove this theorem 
it is helpful to strengthen slightly the requirement that the linearisation is adapted. This strengthening does not alter $X^{s,{\hat{U}}}_{\min+}$ or $X^{s,{\hat{H}}}_{\min+}$ or their quotients $X\env \hU$ and $X\env \hH$. 
The proof  in \cite{BDHK2} (which is then strengthened in \cite{BDHK3}) that, if a suitable version of the condition that semistability coincides with stability is satisfied, the algebras of invariants $\oplus_{m=0}^\infty H^0(X,L^{\otimes cm})^{\hat{U}}$ and $$\oplus_{m=0}^\infty H^0(X,L^{\otimes cm})^{\hat{H}} = (\oplus_{m=0}^\infty H^0(X,L^{\otimes cm})^{\hat{U}})^{R}$$ are  finitely generated (and thus that the enveloping quotients  $X \env \hU = X^{s,{\hat{U}}}_{\min+}/\hat{U}$ and $X \env \hH$ are the associated projective varieties) proceeds by induction on the dimension of $U$ and requires that the linearisation is twisted by a \lq {\em well adapted}' rational character $\chi$. More precisely, it is shown in \cite{BDHK2} that, given a linear action of $\hH$ on $X$ with respect to an ample line bundle $L$, there exists $\epsilon >0$ such that if $\chi$ is a rational character of $\GG_m$ (lifting to  a rational character of $\hat{H}$ with kernel containing $H$) with 
$$\omega_{\min} < \chi < \omega_{\min} + \epsilon,$$
and if a suitable \lq semistability coincides with stability' condition is satisfied, 
then  the algebras of invariants $\oplus_{m=0}^\infty H^0(X,L^{\otimes cm})^{\hat{U}}$ and $\oplus_{m=0}^\infty H^0(X,L^{\otimes cm})^{\hat{H}} $ are  finitely generated, and  the enveloping quotients $X\env \hU$ and $X\env \hH$ are the associated projective varieties with  $X \env \hU$ a geometric quotient of $X^{s,{\hat{U}}}_{\min+}$ and $X \env \hH$ a categorical quotient of $X^{s,{\hat{H}}}_{\min+}$. Here $X \env \hH$ is  the reductive GIT quotient of $X\env \hat{U}$ by the induced action of the reductive group ${\hat{H}}/\hat{U} \cong R$ with respect to the linearisation induced by a sufficiently divisible tensor power of $L$.

\end{rmk}

Applying Theorem \ref{mainthm} with $X$ replaced by $X \times \PP^1$, with respect to the tensor power of the linearisation $L$ (over $X$) with $\mathcal{O}_{\PP^1}(M)$ (over $\PP^1$) for $M>>1$, gives us a  projective variety $(X \times \PP^1) \env \hat{H}$ which is a categorical quotient by $\hat{H}$ of an $\hat{H}$-invariant open subvariety of $X \times \kk$. This open subvariety is the inverse image in $(X \times \PP^1)^{s,{\hat{U}}}_{\min+}$ of the $R$-semistable subset $((X\times \PP^1) \env \hat{U})^{ss,R}$ of $(X\times \PP^1) \env \hat{U} = (X \times \PP^1)^{s,{\hat{U}}}_{\min+}/{\hat{U}}$,
and contains as an open subvariety a geometric quotient by $H$ of an $H$-invariant open subvariety $X^{\hat{s},H}$
 of $X$. 

\begin{rmk}
Here $X^{\hat{s},H}$ can be identified in the obvious way with $X^{\hat{s},H}  \times \{ [1:1]\} $ which is the intersection with $X \times \{ [1:1]\}$ of the inverse image in $(X \times \PP^1)^{\rms,{\hat{U}}}_{\min+} = (X \times \PP^1)^{\ssfg,{\hat{U}}}_{\min+}$ of the $R$-stable subset $((X\times \PP^1) \env \hat{U})^{s,R}$ of 
$$(X\times \PP^1) \env \hat{U} =  ((X^0_{\min} \times \kk^*) \sqcup (X^{s,{\hat{U}}}_{\min+}  \times \{0\}))/{\hat{U}} \cong 
(X^0_{\min}/U) \sqcup (X^{s,{\hat{U}}}_{\min+}  /{\hat{U}}).$$
This geometric quotient $X^{\hat{s},H}/H$ and its projective completion $(X \times \PP^1) \env \hat{H}$ can be described using Hilbert--Mumford-like criteria, by combining the description of 
$(X\times \PP^1) \env \hat{U}$ as the geometric quotient $ (X \times \PP^1)^{s,{\hat{U}}}_{\min+}/{\hat{U}}$ with reductive GIT for the induced linear action of the reductive group $R=H/U$ on $(X\times \PP^1) \env \hat{U}$.
\end{rmk} 

Theorem \ref{mainthm} describes the good case when semistability coincides with stability for the linear action of $\hU$. Theorem \ref{mainthm2} below, which is proved in \cite{BDHK3},  applies to any adapted linear action of  $\hat{H}$, provided that the much weaker condition that the stabiliser in the unipotent radical $U$ of a generic $x \in X$ is trivial. 

\begin{rmk} \label{remark1.11}
In fact this weaker hypothesis can itself be removed. It is shown in \cite{BDHK3} that Theorem \ref{mainthm} is still true when condition ($\mf{C}^*$), that semistability coincides with stability for $\hU$, is replaced with the weaker condition 
\begin{equation}  
\text{$\dim(\stab_{U}(z)) = \min_{x \in X} \, \dim(\stab_{U}(x)) $ for every $z \in Z_{\min}$}. \tag{$\tilde{\mf{C}}^*$}
\end{equation}
Theorem \ref{mainthm2} is then still valid without the hypothesis that the stabiliser in the unipotent radical $U$ of a generic $x \in X$ is trivial, provided that the condition ($\mf{C}^*$) is replaced with  ($\tilde{\mf{C}}^*$) in its statement.
\end{rmk}

Theorem \ref{mainthm2} is a non-reductive analogue of the partial desingularisation construction for reductive GIT described at the end of $\S$1.1.

\begin{thm} \label{mainthm2} 
Let $H$ be a linear algebraic group over $\kk$ with graded unipotent radical $U$
 and let $\hat{H} = H \rtimes \GG_m$ be the extension of $H$ by $\GG_m$ 
which defines the grading.
Suppose that $\hH$  acts linearly on an irreducible projective variety $X$ with respect to an adapted ample linearisation. Suppose also that $\stab_U(x) = \{e\}$ for generic $x \in X$.

Then there is a sequence of blow-ups of $X$ along $\hat{H}$-invariant projective subvarieties (the first of which is the closure in $X$ of the locus where the stabiliser in $U$ has maximal dimension in $X^0_{\min}$)
 resulting in a projective variety $\hat{X}$ with an adapted linear action of $\hat{H}$ (with respect to a power of an ample line bundle given by tensoring the pullback of $L$ with small multiples of the exceptional divisors for the blow-ups) which satisfies the  condition  ($\mf{C}^*$), so that Theorem \ref{mainthm} applies.

Moreover there is a sequence of further blow-ups along $\hat{H}$-invariant projective subvarieties appearing as the closures of $H$-sweeps of connected components of fixed point sets of reductive subgroups of $H$, resulting in a projective variety $\tilde{X}$ satisfying the same conditions as $\hat{X}$ and in addition that the enveloping quotient $\tilde{X}\env \hat{H} $ is the geometric quotient by $\hH$ of  the  $\hH$-invariant open subset $\tilde{X}^{s,\hH}_{\min+}$. Both $\hat{X}\env \hat{H}$ and $\tilde{X}\env \hat{H}$ are projective completions of the geometric quotient by $\hH$ of the $\hH$-invariant open subset $X^{s,\hH}_{\min+}$ of $X$ which can be identified via the blow-down map with the complement in $\tilde{X}^{s,\hH}_{\min+}$ of the exceptional divisors.
\end{thm}

By considering the action of $\hH$ on $\widehat{X \times \PP^1}$ (and similarly on $\widetilde{X \times \PP^1}$)
as above, we obtain a projective variety $\widehat{X \times \PP^1} \env \hat{H}$ which is a categorical quotient by $\hat{H}$ of a $\hat{H}$-invariant open subset of a blow-up of ${X} \times \kk$ and contains as an open subset a geometric quotient of an $H$-invariant open subset $X^{\hat{s},H}$ of $X$ by $H$, where the geometric quotient $X^{\hat{s},H}/H$ and its projective completion $\widehat{X \times \PP^1} \env \hat{H}$  have descriptions in terms of Hilbert--Mumford-like criteria, the explicit blow-up construction used to obtain $\widehat{X \times \PP^1} $ from $X \times \PP^1$ and an analogue of S-equivalence.

\section{Extended, graded and torus-graded linearisations}

Recall from Definition \ref{defnextendedlin} that an extended linearisation $\mathcal{L}$ of an action of a linear algebraic group $H$ on a projective variety $X$ is given by the data:\\
(a) a line bundle $L$ on $X$ (usually assumed to be ample);\\
(b) a semi-direct product $\hat{H} = H \rtimes \GG_m$ of $H$ by $\GG_m$; \\
(c) an extension of the $H$-action on $X$ to $\hat{H}$ and a lift of the $\hat{H}$-action to $L$.

Recall also that given an extended linearisation $\mathcal{L}$ and a rational number $q \in \QQ$, we define the $q$-hat-stable locus $X^{\hat{s},q} = X^{\hat{s},q,\mathcal{L}}$ to be the $H$-{invariant open subvariety}  of $X$ determined by
$$  X^{\hat{s},q} \times  \{ [1:1]\}  = (X \times \PP^1)^{s,\hat{H}} \cap (X \times \{ [1:1]\} )$$
where $\hH$ acts on $\PP^1$ as at 
(\ref{p1action}) above  with its linearisation on $\mathcal{O}_{\PP^1}(1)$ twisted by $q$, and 
$(X \times \PP^1)^{s,\hat{H}}$ is defined 
with respect to the induced linearisation for the $\hH$-action on $X \times \PP^1$ on $L$ tensored with $\mathcal{O}_{\PP^1}(m)$ for $m\! >\!>\! 1$.
We then have a {geometric quotient} 
$X^{\hat{s},q}/H$   
by $H$ which is isomorphic to the open subset $ ( (X \times \PP^1)^{s,\hat{H}} \cap (X \times (\kk \setminus \{ 0 \}) ))/\hH
$ of $(X \times \PP^1)^{s,\hat{H}} / \hat{H}$ for this choice of linearisation. 

\begin{rmk}
Given a linearisation in the classical sense of an action of a linear algebraic group $H$ on a projective variety $X$ with respect to a line bundle $L$, we have a \lq trivial extended linearisation' for which $\hat{H} = H \times \GG_m$ and $\GG_m$ acts trivially on $X$ and on $L$.  If $q \in (0,1)$ then the stable locus for the action of $\GG_m$ on $X \times \PP^1$
with respect to the induced linearisation for the $\hH$-action on $X \times \PP^1$ on $L$ tensored with $\mathcal{O}_{\PP^1}(m)$
 is 
$X \times (\kk \setminus \{ 0 \})$. 
Thus taking $m\! >\!>\! 1$ 
the $q$-hat-stable locus $X^{\hat{s},q}$ coincides with the stable locus defined as in \cite{BDHK} for the action of $H$ on $X$ with the given linearisation. Similarly  if $q \notin [0,1]$ the $q$-hat-stable locus $X^{\hat{s},q}$ is empty for this linearisation.
\end{rmk}

Recall from Definition \ref{defngradedlin}  that a  linear algebraic group with graded unipotent radical is a linear algebraic group $H$ with unipotent radical $U$, equipped with a semi-direct product $\hH = H \rtimes \GG_m$ such that the adjoint action of $\GG_m$ on the Lie algebra of $U$ has only strictly positive weights and the induced conjugation action of $\GG_m$ on $H/U$ is trivial. 
Recall also that a graded linearisation $\mathcal{L}$ of an action of  $H$ on a projective variety $X$ is then an extended linearisation in the sense of Definition \ref{defnextendedlin} for this choice of $\hH$, such that the $\hH$-linearisation is twisted  by a rational character of $\hH$ so that 0 lies in the interior of the lowest bounded chamber for the $\GG_m$ action.  Given a graded linearisation $\mathcal{L}$,  the \lq hat-stable' locus $X^{\hat{s}} = X^{\hat{s},\mathcal{L}}$ is the $0$-hat-stable locus  $X^{\hat{s},0}$ as defined in Definition \ref{defnextendedlin} when $q=0$.

\begin{rmk} \label{naturalgrad}
When $H$ is a linear algebraic group with graded unipotent radical $U$ and $\mathcal{L}$ is a graded linearisation for an action of $H$ on a projective variety $X$ (with respect to an ample line bundle $L$ on $X$), then we can apply Theorems \ref{mainthm} and \ref{mainthm2} to the action of $\hH$ on $X \times \PP^1$. Thus the {geometric quotient} $X^{\hat{s}}/H$ by $H$ has a {projective completion} which
 is a geometric quotient by $\hat{H}$ of an open subset of a $\hat{H}$-equivariant blow-up of $X \times  \PP^1$. 
Furthermore  the {geometric quotient} $X^{\hat{s}}/H$ by $H$ and its {projective completion} can be described using
Hilbert--Mumford-like criteria combined with the explicit blow-up construction. Thus the data of the graded linearisation gives us a  GIT-like quotient with most of the good properties which hold in the reductive case.
\end{rmk}

Now  let $T$ be a torus defined over $\kk$.
Definitions \ref{defnextendedlin}  and \ref{defngradedlin} can be generalised to define $T$-extended linearisations, and $T$-graded linearisations for the actions of linear algebraic groups with $T$-graded unipotent radical.

 \begin{definition} \label{defnextendedlin2} A {\em $T$-extended linearisation} $\mathcal{L}$ of an action of a linear algebraic group $H$ on a projective variety $X$ is given by the data:\\
(a) a line bundle $L$ on $X$ (usually assumed to be ample);\\
(b) a semi-direct product $\hat{H} = H \rtimes T$ of $H$ by $T$; \\
(c) an extension of the $H$-action on $X$ to $\hat{H}$ and a lift of the $\hat{H}$-action to $L$.\\
Given an extended linearisation $\mathcal{L}$ for the action of $H$ on $X$, and a projective toric variety $Y = \overline{Ty_0}$ 
 with an ample linearisation $\mathcal{L}_T$ for the action of $T$ on $Y$, we can define the \lq $(Y,\mathcal{L}_T)$-hat-stable' locus $ X^{\hat{s},Y,\mathcal{L}_T,\mathcal{L}}$ to be the $H$-{invariant open subvariety}  of $X$ determined by
$$  X^{\hat{s},Y,\mathcal{L}_T,\mathcal{L}} \times  \{ y_0\}  = (X \times Y)^{s,\hat{H}} \cap (X \times \{y_0\} )$$
where $(X \times Y)^{s,\hat{H}}$ is defined as in \cite{BDHK} 
with respect to the induced linearisation for the $\hH$-action on $X \times Y$ with respect to the linearisation $\mathcal{L}$ tensored with $\mathcal{L}_T^{\otimes m}$ for $m\! >\!>\! 1$.
\end{definition}

\begin{rmk}
We then have a {geometric quotient} $X^{\hat{s}, Y,\mathcal{L}_T,\mathcal{L}}/H$ by $H$ which is isomorphic to an open subset of $(X \times Y)^{s,\hat{H}} / \hat{H}$ for a choice of linearisation as in Definition \ref{defnextendedlin2}. 
\end{rmk}

\begin{definition} \label{defngradedlin2}
A {\em  linear algebraic group with $T$-graded unipotent radical} is a linear algebraic group $H$ with unipotent radical $U$, equipped with \\
i) a semi-direct product $\hH = H \rtimes T$  such that  the induced conjugation action of $T$ on $H/U$ is trivial, and \\
ii) a non-empty open rational cone $C$ in the Lie algebra of $T$ such that the adjoint action on the Lie algebra of $U$ of any one-parameter subgroup of $T$ whose derivative at the identity lies in $C$ has only strictly positive weights. 

A {\em $T$-graded linearisation} $\mathcal{L}$ of an action of  $H$ on a projective variety $X$ is then a $T$-extended linearisation in the sense of Definition \ref{defnextendedlin2} for this choice of $\hH$, with the $\hH$-linearisation twisted  by a rational character of $\hH$ whose kernel contains $H$, in such a way that 0 lies in the interior of the lowest bounded chamber for some 
 one-parameter subgroup of $T$ whose derivative at the identity lies in the cone $C$.  
\end{definition}

 When $T$ is the one-parameter multiplicative group  $\GG_m$ and $Y=\PP^1$, and $\mathcal{L}_T$ is  the linearisation of the $\GG_m$-action on $\mathcal{O}_{\PP^1}(1)$ given by the representation (\ref{p1action}), 
 then we recover the definitions of extended and graded linearisations given above.

\begin{rmk}
When $H$ is a linear algebraic group with $T$-graded unipotent radical $U$ and $\mathcal{L}$ is a $T$-graded linearisation for an action of $H$ on a projective variety $X$ with respect to an ample line bundle $L$ on $X$, then 
an analogous picture to that of Remark \ref{naturalgrad} holds \cite{ bejk16}. Thus the {geometric quotient} $X^{\hat{s}, Y,\mathcal{L}_T,\mathcal{L}}/H$  has a {projective completion} which
 is a geometric quotient by $\hat{H}$ of an open subset of a $\hat{H}$-equivariant blow-up of the product of $X$ with the toric variety $ Y$. 
Furthermore  the {geometric quotient} $X^{\hat{s}, Y,\mathcal{L}_T,\mathcal{L}}/H$  and its {projective completion} can be described using
Hilbert--Mumford-like criteria combined with the geometry of the toric variety $T$, the rational cone $C$ and the  blow-up construction. 
\end{rmk}

\section{Applications} 

In this section we will describe some linear actions of non-reductive groups where GIT for suitable graded linearisations, obtained as in Remark \ref{naturalgrad}, behaves better than GIT for classical linearisations. 

\begin{ex}
The first of these are the famous Nagata counterexamples to Hilbert's 14th problem \cite{Nagata}, which provide examples of linear actions of unipotent groups $U$ on  projective space such that the corresponding algebras of $U$-invariants are not finitely generated. In these examples the linear action extends to a linear action of an extension $\hat{U} = U \rtimes \GG_m$ by $\GG_m$
such that the action of $\GG_m$ by conjugation on the Lie algebra of $U$ has all its weights strictly positive, and $\stab_U(x) = \{e\}$ for generic $x$, so Theorem \ref{mainthm2} applies, and the quotient $(X \times \PP^1)\env \hU$ gives us a projective completion of a geometric quotient by $U$ of a $U$-invariant open subset of the projective space $X$ which can be determined by Hilbert--Mumford-like criteria. We can regard this as the GIT quotient of the projective space by the graded unipotent group $U$ with respect to the induced graded linearisation. Note that the $\hU$-invariants on $X \times \PP^1$ restrict to $U$-invariants on $X$, so the grading is picking out for us a finitely generated subalgebra of the algebra of $U$-invariants, and thus a tractable GIT quotient.
\end{ex}

\begin{ex}
Recall that the automorphism group of the weighted projective plane $\PP(1,1,2) = (\kk^3 \setminus \{ 0 \})/ \GG_m$, for $\GG_m$ acting linearly on $\kk^3$ with weights $1,1,2$,  is given by 
$$\mbox{Aut}(\PP(1,1,2)) \cong R \ltimes U$$
where $R \cong \GL(2)$ is reductive and 
$U \cong (\GG_a)^3$  is unipotent, with elements  $(\lambda,\mu,\nu) \in (\kk)^3$ acting on $\PP(1,1,2)$ via 
$$[x,y,z]  \mapsto [x,y,z+\lambda x^2 + \mu xy + \nu y^2].$$
The central one-parameter subgroup $\GG_m$ of $R \cong \GL(2)$ acts on the Lie algebra of 
$U$ with all positive weights, and the
associated semi-direct product
$$\hat{U} = U \rtimes \GG_m$$
can be identified with a subgroup of $\mbox{Aut}(\PP(1,1,2))$. Thus any ample linearisation for an action of $\mbox{Aut}(\PP(1,1,2))$ on a projective variety $X$ becomes a graded linearisation in a natural way. It therefore follows from Theorem \ref{mainthm2} that whenever $H=\mbox{Aut}(\PP(1,1,2))$ acts linearly on a projective
variety $X$
and $\stab_U(x) = \{e\}$ for generic $x \in X$, then  there is a geometric quotient by $H$ of an open subset of $X$ described by Hilbert--Mumford-like criteria, with a projective completion  which is a categorical quotient  of an open subset  of an $H$-equivariant blow-up $\tilde{X}$ of $X$.

Indeed the same is true for the automorphism group of { any} complete simplicial toric variety.
For it was observed in \cite{BDHK2} using the description in \cite{cox} that
the automorphism group $H$ of
any complete simplicial toric variety is a linear algebraic group with a graded unipotent radical $U$; there is a grading defined by a one parameter subgroup $\GG_m$ of $H$ acting by conjugation on the Lie algebra of $U$ with all weights strictly positive, and inducing a central one-parameter subgroup of $R=H/U$. Thus Theorems \ref{mainthm} and \ref{mainthm2} (and if necessary Remark \ref{remark1.11}) can be applied.
\end{ex}

\begin{ex}
\label{subsec:reparametrisation groups}

Suppose now that $\kk = \CC$ and consider $k$-jets at 0 of holomorphic maps from $\CC^p$ to a complex manifold $Y$  for any $k, p \geq 1$. It was observed in \cite{BK15} that 
 the group $\GG_{(k,p)}$ of $k$-jets of holomorphic reparametrisations of $(\CC^p,0)$  has 
a graded unipotent radical $\UU_{(k,p)}$ such that the grading is defined by a one-parameter subgroup of $\GG_{(k,p)}$ acting by conjugation on the Lie algebra of $\UU_{(k,p)}$ with all weights strictly positive, and inducing a central one-parameter subgroup of the reductive group $\GG_{(k,p)}/\UU_{(k,p)}$. 
So  Theorems \ref{mainthm} and \ref{mainthm2}, with  Remark \ref{remark1.11}, can be applied 
to any linear action of the
reparametrisation group  $\GG_{(k,p)}$.

\end{ex}

\begin{ex}
\label{subsec:unstable strata}

Finally  let $G$ be a reductive group over an algebraically closed field $\kk$ of characteristic zero, acting linearly on a projective variety $X$ with respect to an ample line bundle $L$. Associated to this linear $G$-action and an invariant inner product on the Lie algebra of $G$, there is a stratification 
$$ X = \bigsqcup_{\beta \in \mathcal{B}} S_\beta$$ of $X$ by locally closed subvarieties $S_\beta$, 
indexed by a partially ordered finite subset $\mathcal{B}$ of a positive Weyl chamber for the reductive group $G$,  such that 

 (i) $S_0 = X^{ss}$, 

\noindent and for each $\beta \in \mathcal{B}$

 (ii) the closure of $S_\beta$ is contained in $\bigcup_{\gamma \geqslant \beta} S_\gamma$, and

 (iii) $S_\beta \cong G \times_{P_\beta} Y_\beta^{ss}$

\noindent where
$P_\beta$ is a parabolic subgroup of $G$ acting on  a projective subvariety $\overline{Y}_\beta$ of $X$ with an open subset $Y_\beta^{ss}$ which is determined by the action of the Levi subgroup $L_\beta$ of $P_\beta$ with respect to a suitably twisted linearisation \cite{Hess,K}.

Here the original linearisation for the action of $G$ on $L \to X$ is restricted to the action of the parabolic subgroup $P_\beta$ over $\overline{Y}_\beta$, and then twisted by the rational character $\beta$ of $P_\beta$ which is almost adapted for a central one-parameter subgroup of the Levi subgroup  $L_\beta$ acting with all weights strictly positive on the Lie algebra of the unipotent radical of $P_\beta$. So $P_\beta$ is a linear algebraic group with graded unipotent radical; indeed, its unipotent radical is graded by the torus which is the centre of a Levi subgroup of $P_\beta$.
Thus to construct a quotient by $G$ of (an open subset of) an unstable stratum $S_\beta$, we can study the linear action on $\overline{Y}_\beta$ of the parabolic subgroup $P_\beta$, and apply   Theorems \ref{mainthm} and \ref{mainthm2}, with Remark \ref{remark1.11}.

In this situation $Y_{\beta} $ is equal to $(\overline{Y}_\beta)^0_{\min}$ (in the notation introduced in $\S$1.3, immediately before Theorem \ref{mainthm}), and we have a retraction $p_{\beta}: Y_{\beta} \to Z_{\beta}$ where 
$$ p_{\beta}(y) =  \lim_{t \to 0, \,\, t \in \GG_m} t \cdot y
$$
for $y \in Y_{\beta}$ and $Z_{\beta}$ plays the role of $Z_{\min}$ in $\S$1.3. Since the rational character $\beta$ of $P_{\beta}$ is borderline adapted, not adapted, we have 
$$Y_{\beta}^{ss} = p_{\beta}^{-1}(Z_{\beta}^{ss,L_{\beta}})$$
and the reductive GIT quotient $Z_{\beta} /\!/ L_{\beta}$ is the categorical quotient of $Y_{\beta}^{ss}$ by $P_{\beta}$. However this is certainly not a geometric quotient (because the closure of every $P_{\beta}$-orbit in  $Y_{\beta}^{ss}$ meets  $Z_{\beta}^{ss}$).  $Z_{\beta} /\!/ L_{\beta}$ is also the categorical quotient of $P_{\beta}Z_{\beta}^{ss,L_{\beta}} \subseteq Y_{\beta}^{ss}$ by $P_{\beta}$, and 
$$P_{\beta}Z_{\beta}^{ss,L_{\beta}}  = U_{\beta}Z_{\beta}^{ss,L_{\beta}}  $$
where $U_{\beta}$ is the unipotent radical of $P_{\beta}$.

On the other hand if we modify the linear action of $P_{\beta}$ on  $\overline{Y}_\beta$ by an adapted rational character, given by $(1+ \delta)\beta$ for $0 < \delta <\!< 1$, rather than by $\beta$, then in the notation of $\S$1.3 we have
$$(\overline{Y}_\beta)^{s,P_{\beta}}_{\min +} = Y_{\beta}^{ss} \setminus U_{\beta}Z_{\beta}^{ss,L_{\beta}}.$$
Thus, at least after applying suitable blow-ups and removing any resulting unstable strata, we can use  Theorems \ref{mainthm} and \ref{mainthm2} (and if necessary Remark \ref{remark1.11}) to construct a categorical quotient by $P_\beta$ of a $P_\beta$-invariant open subvariety of (a blow-up of) $Y_{\beta}^{ss} \setminus U_{\beta}Z_{\beta}^{ss,L_{\beta}}$ which will fibre over $Z_{\beta} /\!/ L_{\beta}$.

There are many moduli spaces (of \lq stable' or \lq semistable' objects) in algebraic geometry which can be constructed as GIT quotients of reductive group actions on projective varieties with respect to ample linearisations. We can hope to use the construction just described to construct corresponding moduli spaces of unstable objects.
In particular we can consider moduli spaces of sheaves of fixed Harder--Narasimhan type over a nonsingular projective variety $W$ (cf.  \cite{hok12}). 
There are  well known constructions going back to Simpson \cite{s94}
of the  moduli spaces of semistable pure sheaves on $W$ of fixed  Hilbert polynomial as GIT
quotients of linear actions of suitable  special linear groups $G$ on  schemes $Q$ (closely related to 
quot-schemes) which are $G$-equivariantly embedded in projective spaces.
 These constructions can
be chosen so that elements of $Q$ which parametrise sheaves of a fixed Harder--Narasimhan
type  form a stratum in the stratification of $Q$ associated to the linear action of $G$ (at least modulo
taking connected components of strata) \cite{hok12}. The associated linear actions of  parabolic subgroups of these special linear groups $G$ have natural graded linearisations, and so can be used to construct and study (projective completions of) moduli spaces of  sheaves
of fixed Harder--Narasimhan type over $W$ \cite{behjk16}.
The simplest non-trivial case is that of unstable vector bundles of rank 2 and fixed Harder--Narasimhan type over a nonsingular projective curve $W$ (cf. \cite{bramn09}).

The other classical moduli spaces constructed as GIT quotients are the moduli spaces of stable curves of fixed genus $g \geq 2$. Here too it is possible to construct quotients of suitable \lq unstable strata', giving us moduli spaces of unstable curves when appropriate discrete invariants of the singularities are fixed; these play a role analogous to that of Harder--Narasimhan types \cite{jj}. Moreover in the case of curves of genus $g \geq 2$ the condition that semistability coincides with stability for the extension of the unipotent radical is always satisfied, so we can apply Theorem \ref{mainthm} without having to resort to the partial desingularisation construction of Theorem \ref{mainthm2}. 
\end{ex}

\end{document}